\documentclass[oneside]{amsart}
\usepackage[mathscr]{eucal}
\usepackage{amssymb, amsmath,array, amscd}

\theoremstyle{change}

\newtheorem{THM}{Theorem}

\newtheorem{prop}{Proposition}[section]

\newtheorem{lemma}{Lemma}[section]

\newtheorem{definition}{Definition}[section]

\newtheorem{remark}{Remark}

\newtheorem{cor}{Corollary}[section]

\newtheorem{example}{Example}[section]

\begin{document}

\title[Fibrations, Divisors and Transcendental Leaves]
{Fibrations, Divisors and Transcendental Leaves}

\author{Jorge Vit\'orio Pereira  \\ \, \\ (with an appendix by  Laurent Meersseman) }

\begin{abstract}
We will use flat divisors, and  canonically associated singular
holomorphic foliations, to investigate some of the geometry of
compact complex manifolds. The paper is mainly concerned with
three distinct problems: the existence of fibrations, the topology
of smooth hypersurfaces and the topological closure of
transcendental leaves of foliations.
\end{abstract}

\thanks{The author was partially supported by PROFIX/Cnpq and  CNRS}



\maketitle

\tableofcontents

\section{Introduction and Statement of Results}

Let $M$ be a compact complex manifold and ${\rm Div(M)}$ be its group of divisors, i.e., the free abelian group
generated be the irreducible hypersurfaces of $M$. We will denote by ${\rm Div}_{S^1}(M)$ the subgroup of ${\rm Div}(M)$
formed by $S^1$ -flat divisors, i.e., divisors $D$ whose associated line bundle $\mathcal O_M(D)$ admits a hermitian
flat  connection. We will denote by $\Gamma(M)$ the quotient of the group of rational divisors by the group of rational
$S^1$-flat divisors, i.e.,
\[
\Gamma(M) = \frac{{\rm Div(M)}\otimes \mathbb Q}{{\rm Div}_{S^1}(M) \otimes \mathbb Q }.
\]
If $M$ is a projective manifold then
\[
  \Gamma(M) = NS_{\mathbb Q}(M) \, ,
\]
i.e., $\Gamma(M)$ can be identified
with the rational Neron-Severi group of $M$, see section \ref{S:Gamma}, and in particular is finite dimensional.
Our first result says that this is always the case for compact complex manifolds, i.e.,

\begin{THM}\label{T:finito}
Let $M$ be a compact complex manifold then $\dim_{\mathbb Q} \Gamma(M) < \infty$.
\end{THM}

Most of our results involve {\it the $\Gamma$-class of a divisor},
i.e., the image of the divisor under the natural homomorphism
\[
   {\rm Div(M)} \to \Gamma(M) \, .
\]

Note that the $\Gamma$-class of a divisor completely determines
its rational-Chern class or, equivalently,  its rational homology
class. For Kaehler varieties the converse also holds, i.e., the
$\Gamma$-class of a divisor is completely determined by its
rational Chern-class. In fact in this case $\Gamma(M)$ is
isomorphic to the image of the group of divisors under the
rational Chern class map.

>From some basic properties of flat divisors and  Theorem
\ref{T:finito} we  deduce  a characterization of compact complex
manifolds which fibers over a projective curve.

\begin{THM}\label{T:fibracao}
Let $M$ be a compact complex  manifold with  $\dim_{\mathbb Q} \Gamma(M) + 2$  pairwise disjoint hypersurfaces $H_i$. Then there exists a
holomorphic map $\rho:M\to C$, $C$ a smooth algebraic curve, with connected fibers such that every $H_i$ is a component
of a fiber of $\rho$.
\end{THM}

This result was inspired by an unpublished   result due to A. Vistoli stating that a smooth  projective variety,
over an arbitrary field,  with an infinite number of disjoint hypersurfaces
fibers over a curve.  In the case of compact Kaehler varieties Vistoli's proof has been adapted
by M. Sebastiani in \cite{Sebastiani}.

A related result for projective varieties has been obtained by B. Totaro, see Theorem 2.1 of \cite{Totaro}. The proof of Theorem \ref{T:fibracao}
also prove a generalization of Theorem 2.1 of \cite{Totaro}  to arbitrary complex manifolds.

\begin{THM}\label{T:Totaro2}
Let $M$ be a compact complex variety.
Let $D_1, D_2, \ldots, D_r$, $r\ge 3$, be connected effective divisors which are pairwise disjoint and whose
$\Gamma$-classes lie in a line of  $\Gamma(M)$. Then there exists a map $\rho: M\to C$ with connected
fibers to a smooth curve $C$ such that $D_1, \ldots, D_r$ which maps the divisors $D_i$ to points.
\end{THM}

As we will see some real foliations of  codimension one  are
naturally associated to flat divisors. The pertubation of these
foliations allow us to prove  the next result, 

\begin{THM}\label{T:Totaro}
Let $M$ be a compact complex manifold.
If $H_1$ and $H_2$ are smooth connected disjoint  hypersurfaces  such that $[H_1]$ and $[H_2]$ lie in the same line of $\Gamma(M)$  then there exists
an \'etale $\mathbb Z/n$-covering $\tilde{D_1}$ of $D_1$ and an \'etale $\mathbb Z/m$-covering $\tilde{D_2}$ of $D_2$ which are
diffeomorphic where $m$ and $n$ are positive integers satisfying $m[H_1] = n[H_2]$.
\end{THM}

 It has to be noted that when $M$ and $H_1$ and $H_2$ are smooth connected hypersurfaces 
 and the Picard variety of $M$ is isogeneous to a product of elliptic curves then it is shown in \cite{Totaro} that there exists
 finite and cyclic \'etale coverings of $H_1$ and $H_2$ with the same pro-l homotopy type. Theorem \ref{T:Totaro} is a generalization of this
 result and gives a positive answer to a conjecture made by B. Totaro in the above mentioned paper.

In analogy with the projective case we will say that a leaf $L$ of a codimension one holomorphic foliation of a compact complex manifold $M$ is
transcendental if it is not contained in any compact complex hypersurface.  Our last main result is the following

\begin{THM}\label{T:FolhasGeral}
Let $M$ be a compact complex manifold, $\mathcal F$ a codimension
one holomorphic foliation of $M$ and $L$  a transcendental leaf of
$\mathcal F$. Denote by $\mathcal H$ the  set  of compact complex
irreducible hypersurfaces of $M$ which do not intersect the
topological closure of $L$. Then the following assertion holds
\begin{enumerate}
\item In general the cardinality of $\mathcal H$ is at most $\dim_{\mathbb Q} \Gamma(M) +
1$ and when is equal to $\dim_{\mathbb Q} \Gamma(M) + 1$ then
$\mathcal F $ is given by a closed meromorphic $1$-form;
\item If $M$ is projective then  the cardinality of $\mathcal H$ is at most $\dim_{\mathbb Q} \Gamma(M) $
 and when is equal to $\dim_{\mathbb Q} \Gamma(M) $ then
 $h^1(M,\mathcal O_M)\neq 0$ and $\mathcal F$ is given by a closed meromorphic
 $1$-form. In particular for projective manifolds without global
 holomorphic $1$-forms we have that the cardinality of $\mathcal
 H$ is at most $\dim_{\mathbb Q} \Gamma(M) - 1$.
\end{enumerate}
\end{THM}

We remark that we cannot replace in the statement of Theorem
\ref{T:FolhasGeral} item (1) the group $\Gamma(M)$ by the
Neron-Severi group of $M$, see Section \ref{S:Folhas}. In an
appendix to this paper L. Meersseman constructs a complex manifold
of dimension $5$ showing that in Theorems \ref{T:fibracao},
\ref{T:Totaro2} and \ref{T:Totaro} it is also not possible to
replace $\Gamma(M)$ by the Neron-Severi group of $M$.

\vskip.3cm

\noindent{\bf Acknowledgements:} I am largely indebted to Steven
Kleiman who  brought to my attention   the results  of A. Vistoli
and B. Totaro. I am also indebted to  Marcos Sebastiani who
explained to me  some properties of the cohomology of compact
complex manifolds and to  Burt Totaro whose comments on a
preliminary version of the present paper allowed me to clarify the
arguments used to prove Theorem \ref{T:Totaro}. Conversations with
Marco Brunella and Frank Loray
 at the early stages of this work played a crucial role on its
further development.

 \tableofcontents

\section{Flat line bundles over complex manifolds}\label{S:flat}

If  $M$ is  a complex manifold then the set of isomorphism classes of holomorphic line-bundles is identified with
$H^1(M,{\mathcal O^*_M})$  as follows: let $\mathcal L$ be a holomorphic line bundle,
$\mathcal U = \{ U_i \}$ a sufficiently fine   open covering of $M$ and  $\phi_i: U_i \to \mathcal L$ be nowhere vanishing local
holomorphic sections of $\mathcal L$; over $U_i \cap U_j$ we have  the transition functions $\phi_{ij}=\phi_i \cdot \phi_{j}^{-1}: U_i\cap
U_j \to \mathbb C^*$ satisfying the identities $\phi_{ij}\cdot \phi_{ji} = 1$ and $\phi_{ij}\cdot \phi_{jk}\cdot \phi_{ki} = 1$.
Therefore the collection $\{ \phi_{ij} \}$ determines and element of $H^1(M,{\mathcal O^*_M})$. It can be  verified that
this element does not depend on the choices made above.

Denote by $\mathbb C^*$, resp. $S^1$, the constant sheaf over $M$ of invertible complex numbers, resp. complex numbers of modulus $1$.

\begin{definition}\rm
A line bundle $\mathcal L \in H^1(M,{\mathcal O_M}^*)$ is $\mathbb C^*$-flat, resp. $S^1$-flat, if $\mathcal L$ belongs to
the image of the morphism $H^1(M,{\mathbb C}^*) \to H^1(M,{\mathcal O^*_M})$, resp. $H^1(M,{S}^1) \to H^1(M,{\mathcal O^*_M})$,
induced by the natural inclusions.
\end{definition}

Concretely a line bundle bundle $\mathcal L$ is $\mathbb C^*$-flat, resp. $S^1$-flat, if it admits a system of local holomorphic sections
whose transition functions are locally constant holomorphic functions, resp. locally constant holomorphic functions of modulus $1$.

We recall that  the {\it Chern class} of a line bundle $\mathcal L$, denoted by $c(\mathcal L)$, is the image
of $\mathcal L$ under the map $H^1(M,{\mathcal O^*_M}) \to H^2(M,\mathbb Z)$ induced by the exponential sequence
\[
  0 \to \mathbb Z \to \mathcal O_M \to {\mathcal O^*_M} \to 0 \, .
\]
The {\it real Chern class} of $\mathcal L$, denoted by $c_{\mathbb R}(\mathcal L)$, is the image of $c(\mathcal L)$ under the natural
map $H^2(M,\mathbb Z) \to H^2(M,\mathbb R)$.

The relations between flat line bundles and line bundles with zero real Chern class are presented in the next proposition.
Its content is quite standard but due to a lack of references we will sketch its proof.

\begin{prop}\label{P:flat}
If  $\alpha:H^1(M,\mathbb C) \to H^1(M,\mathcal O_M)$ and $\beta:H^1(M,\mathbb R) \to H^1(M,\mathcal O_M)$ are  the morphisms
induced by the natural inclusions then the following assertions hold:
\begin{enumerate}
\item if $\mathcal L$ is $\mathbb C^*$-flat then $c_{\mathbb R}(\mathcal L) = 0$;
\item if $\alpha$ is surjective then $\mathcal L$ is $\mathbb C^*$-flat
if, and only if, $c_{\mathbb R}(\mathcal  L)=0$;
\item if $\beta$ is surjective then $\mathcal L$ is $S^1$-flat
if, and only if, $c_{\mathbb R}(\mathcal  L)=0$;
\item if the image of $\alpha$ is equal to the image of $\beta$
then $\mathcal L$ is $\mathbb C^*$-flat if, and only if, $\mathcal L$ is $S^1$-flat.
\item if $M$ is compact then $\beta$ is  injective and consequently $\beta$ is an isomorphism if and only if $h^1(M,\mathcal O_M)= 2 h^1(M,\mathbb R)$.
\item if $M$ is compact then  the morphism $H^1(M, S^1) \to H^1(M,\mathcal O_M^*)$, induced by the natural inclusion, is  injective.
\end{enumerate}
\end{prop}
\begin{proof} Consider the commutative diagram of sheaves of abelian groups over $M$
\[
\begin{CD}
0 @>>> \mathbb Z    @>>> \mathbb R @>>> S^1 @>>> 0 \\
 && @VVV  @VVV    @VVV \\
0 @>>> \mathbb Z    @>>> \mathcal O_M @>>> \mathcal O_M^* @>>> 0 \\
 && @AAA  @AAA    @AAA \\
0 @>>> \mathbb Z    @>>> \mathbb C @>>> \mathbb C^* @>>> 0 \\
\end{CD}
\]
with exact rows. From it we obtain the commutative diagram
\[
\begin{CD}
 H^1(M,\mathbb R) @>>> H^1(M,S^1) @>>> H^2(M,\mathbb Z) @>>> H^2(M,\mathbb R) \\
 @V\beta VV  @VVV    @VVV \\
 H^1(M,\mathcal O_M) @>>> H^1(M,\mathcal O_M^*) @>>>  H^2(M,\mathbb Z) \\
  @A\alpha AA  @AAA    @AAA \\
 H^1(M,\mathbb C) @>>> H^1(M,\mathbb C^*) @>>>  H^2(M,\mathbb Z) @>>> H^2(M,\mathbb C)\\
\end{CD}
\]
with exact rows.

The proof of the proposition  will be a standard chasing on the diagram above.

If $\mathcal L$ is a $\mathbb C^*$-flat line-bundle then follows from the exactness of the bottom row of the diagram above that
$c(\mathcal L) \in \ker \{H^2(M,\mathbb Z) \to  H^2(M,\mathbb C) \}$. Since $\ker \{H^2(M,\mathbb Z) \to  H^2(M,\mathbb C) \}=\ker \{H^2(M,\mathbb Z) \to  H^2(M,\mathbb R) \}$
we have that $c_{\mathbb R}(\mathcal L)=0$. This proves assertion (1).

If $\mathcal L$ is a line bundle such that $c_{\mathbb R}(\mathcal L)=0$ then we infer from the diagram that there exists
$\theta \in H^1(M,\mathbb C^*)$ such that $c(\mathcal L \otimes \theta) = 0 \in H^2(M,\mathbb Z)$,
in particular $\mathcal L \otimes \theta \in {\rm Im}\{H^1(M,\mathcal O_M) \to H^1(M,\mathcal O_M^*) \} $.  If $\alpha$ is surjective
it follows that $\mathcal L  \in {\rm Im}\{H^1(M,\mathbb C^*) \to H^1(M,\mathcal O_M^*) \}$ proving that $\mathcal L$ is $\mathbb C^*$-flat.
This proves assertion (2).

Assertions (3) and (4) follows from completely analogous arguments and assertions (5) and (6)  follows from the fact that pluriharmonic functions
over compact complex manifolds are constant.
\end{proof}
 As the proposition above suggests  we do not have in general the equivalence
between  zero real Chern class, $\mathbb C^*$-flat and $S^1$-flat.

\begin{example}\label{E:Hopf}\rm
If $M$ is the quotient of $\mathbb C^2 \setminus \{ 0 \}$ by $(z,w) \mapsto (\lambda_1\cdot z, \lambda_2 \cdot w)$
with $0<|\lambda_1|\le |\lambda_2| < 1$ then $H^1(M,\mathcal O^*_M)=H^1(M,\mathbb C^*)= \mathbb C^*$ and $H^1(M, S^1)= S^1$,
see \cite[pg. 172]{BPV}. In particular every line bundle over $M$ is $\mathbb C^*$-flat and there exists line bundles over $M$ which are not
$S^1$-flat.
\end{example}

\begin{example} \rm
 There exist  complex manifolds $M$, diffeomorphic to  $S^3 \times S^3$, such that $H^2(M,\mathbb Z)=H^1(M,\mathbb C) = 0$ and
$H^1(M,\mathcal O_M) \neq 0 $. Over these manifolds every line bundle has zero real Chern class and a line bundle is $\mathbb C^*$-flat if,
and only if, it is the trivial line bundle.
\end{example}

If $M$ is a compact complex surface then $\alpha$ is always surjective, see \cite[page 117]{BPV}. Therefore it follows
from item (3) of proposition \ref{P:flat} the next

\begin{cor}
On compact complex surfaces a line bundle is $\mathbb C^*$-flat if, and only if, it has zero real Chern class.
\end{cor}

If $M$ is a compact complex Kaehler manifold then it follows from Hodge Theory that $\beta$ is always surjective. In particular we have the following

\begin{cor}
On a compact complex Kaehler manifold a line bundle is $S^1$-flat if, and only if, it has zero real Chern class.
\end{cor}

\section{Flat divisors on   complex manifolds }\label{S:flatdiv}

A divisor $D$ on a complex manifold $M$ is a formal sum
\[
  D = \sum_{i=1}^{\infty} d_i D_i \, ,
\]
where $d_i \in \mathbb Z$ and $\{D_i\}_{ i \in \mathbb N}$ is a locally finite sequence of irreducible hypersurfaces of $M$.
If $\mathcal U = \{ U_i \}$ a sufficiently fine open covering of $M$ then  given a divisor $D$ we can associate a collection of
meromorphic function $f_i: U_i \dashrightarrow \mathbb P^1$ such that the restriction of $D$ to $U_i$  coincides
with the divisor $(f_i)_0 - (f_i)_{\infty}$. Since the functions $f_i$ are unique up to multiplication by a nowhere vanishing holomorphic function
over  $U_i$ we can identify the (additive) group of divisors on $M$, denoted by $Div(M)$, with the (multiplicative) group
$H^0(M,\mathcal M_M^* / \mathcal O_M^*)$, where $\mathcal M_M^*$ denotes the sheaf of invertible meromorphic functions over $M$.

Looking at the  long exact sequence in cohomology associated to the short exact sequence of abelian groups
\[
  0 \to \mathcal O_M^* \to \mathcal M_M^* \to \frac{\mathcal M_M^* }{\mathcal O_M^* } \to 0 \, ,
\]
we obtain a map  $ H^0(M,\mathcal M_M^* / \mathcal O_M^*) \to H^1(M,\mathcal O_M^*)$, i.e., we obtain a map from the group
of divisors over $M$ to the group of isomorphism classes of line bundles over $M$. As usual we will denote the image of a
divisor $D$ by $\mathcal O_M(D)$.

\begin{definition}\rm
We will say that a  divisor $D \in Div(M)$ is $\mathbb C^*$-flat, resp. $S^1$-flat, if  $\mathcal O_M(D)$ is a $\mathbb C^*$-flat, resp. $S^1$-flat, line bundle.
\end{definition}

A divisor $D$ on a complex manifold $M$ is said to be {\it linearly equivalent to zero} if $\mathcal O_M(D) \equiv \mathcal O_M$. Equivalently,
there exits  a meromorphic function $f:M\dashrightarrow \mathbb P^1$  such that
$D = (f)_0 - (f)_{\infty}$. In particular $D$ is $\mathbb C^*$-flat.
Note that if we take $\omega$  the rational $1$-form over $\mathbb P^1$ with simple poles at zero and infinity
then $f^* \omega$ will be a closed meromorphic one-form over $M$
with simple poles along the support of $D$ and holomorphic on the complement. A similar property holds for a general $\mathbb C^*$-flat divisor
as the next proposition shows.

\begin{prop}\label{P:Cflat}
If  $D\neq 0$ is a $\mathbb C^*$-flat line bundle over a complex manifold $M$ then there exists a closed meromorphic one-form $\omega$
with simple poles along the support of $D$ and holomorphic on the complement.
\end{prop}
\begin{proof}
If  $\{ U_i \}$ is sufficiently fine covering of $M$ and $f_i:U_i \dashrightarrow P^1_{\mathbb C}$ are meromorphic  maps  such that
$(f_i)_0 - (f_i)_{\infty} = D_{|U_i}$  then an explicit description of $\mathcal O_M(D)$ is given by the transition functions $f_{ij}=f_i \cdot {f_j}^{-1}:U_i \cap U_j \to
\mathbb C^*$. By hypothesis the   cocycle $\{ f_{ij} \}$  is cohomologous to a locally constant cocycle $\{ t_{ij} \}$. Concretely there exists
holomorphic functions $t_i:U_i \to \mathbb C^*$  such that $f_{ij}\cdot t_{ij}^{-1} = t_i \cdot t_j^{-1}$.
Therefore the  meromorphic   functions
$F_i: U_i \dashrightarrow \mathbb P^1_{\mathbb C}$, $ F_i= \frac{f_i}{t_i}$ satisfies
\begin{enumerate}
\item $(F_i)_0 - (F_i)_{\infty} = D_{|U_i}$
\item $F_i = t_{ij}\cdot F_j$.
\end{enumerate}
In particular, over $U_i\cap U_j$, we obtain   the equality
\[
   \frac{dF_i}{F_i} = \frac{dF_j}{F_j} \, ,
\]
which implies that there exists a closed meromorphic  $1$-form $\omega$ such that
\[
   \omega_{|U_i}=\frac{dF_i}{F_i}  \, ,
\]
for every open set $U_i \in \mathcal U$.
\end{proof}

It is important to note that in general the $1$-form $\omega$ constructed above is not unique. In fact for two distinct choices of flat local
equations for $D$ we obtain two meromorphic $1$-forms with simple poles along $D$ differing by a global closed holomorphic $1$-form $\eta$.

Reciprocally if
$\{F_i: U_i \dashrightarrow \mathbb P^1_{\mathbb C}\}$ is a collection of flat local  equations for $D$ as in the proof of proposition \ref{P:Cflat} and $\eta$ is
a closed holomorphic $1$-form over $M$ then choosing arbitrary branches of $H_i=\exp \int \eta$ over each open set $U_i$ we obtain a new collection of flat local equations
$G_i= F_i/H_i$. Moreover we have that over each open set $U_i$ the equality
\[
 \frac{dG_i}{G_i} = \frac{dF_i}{F_i} - \eta \, .
\]

In other words, if $\Omega^1_{M,{\rm closed}}$ denotes the sheaf of closed holomorphic $1$-forms over $M$ then from the short exact sequence
\[
\begin{CD}
0 @>>> \mathbb C^*     @>>> \mathcal  O_M^* @>d \log >> \Omega^1_{M,{\rm closed}} @>>> 0 \\
\end{CD}
\]
one deduces that the kernel of the map $H^1(M,\mathbb C^*)\to H^1(M,\mathcal O_M ^*)$ coincides with the image of the map
$H^0(M,\Omega^1_{M,{\rm closed}}) \to H^1(M,\mathbb C^*)$.

\begin{example}\label{E:Hopf2}\rm
Let  $M$ be, as in example \ref{E:Hopf},  the quotient of $\mathbb C^2 \setminus \{ 0 \}$ by $(z,w) \mapsto (\lambda_1\cdot z, \lambda_2 \cdot w)$
with $0<|\lambda_1|\le |\lambda_2| < 1$. Suppose further that $\lambda_1^k \neq \lambda_2^l$ for every $(k,l) \in \mathbb Z^2 \setminus \{(0,0)\}$.
Under this hypothesis  $M$ contains just two irreducible curves; they are elliptic curves corresponding to the quotient of the axis by the contraction
above, see \cite[pg. 173]{BPV}. Since  every line bundle over $M$ is $\mathbb C^*$-flat the same is true for divisors. Moreover there is no divisor $D$, $D\neq 0$, on $M$
linearly equivalent to zero.
\end{example}

If we restrict to $S^1$-flat divisors we can refine proposition \ref{P:Cflat} to obtain the stronger

\begin{prop}\label{P:basico}
If  $D\neq 0$ is a $S^1$-flat line bundle over a complex manifold $M$ then there exists a closed meromorphic one-form $\omega_D$
with simple poles along the support of $D$ and holomorphic on the complement. Moreover,
\begin{enumerate}
\item If $H$ is a compact hypersurface which does
not intersect the support of $D$ then $i^* \omega_D = 0$, where $i:H\to M$ denotes the natural inclusion
\item If $M$ is compact the meromorphic $1$-form  $\omega_D$ is unique up
to a multiplicative constant.
\end{enumerate}
\end{prop}
\begin{proof}
As in the proof of proposition \ref{P:Cflat}, and using the same notation, we can construct   meromorphic   functions
$F_i: U_i \dashrightarrow \mathbb P^1_{\mathbb C}$  satisfying
\begin{enumerate}
\item $(F_i)_0 - (F_i)_{\infty} = D_{|U_i}$;
\item $F_i = t_{ij}\cdot F_j$;
\end{enumerate}
where  $t_{ij}$ are locally constant functions. The difference is that now the functions $t_{ij}$ have modulus $1$.

We define $\omega_D$ by the relations
\[
   {\omega_D}_{|U_i}=\frac{dF_i}{F_i}  \, ,
\]
for every open set $U_i \in \mathcal U$.

If $H$  is a  hypersurface disjoint
from the support of $D$ we have that the line bundle $\mathcal O_M(D)$  is trivial
when restricted to $H$.  Moreover if $H$ is compact then the map $H^1(H,S^1)\to H^1(H,\mathcal O_H ^*)$ is injective, see
item (6) of proposition \ref{P:flat},
and it follows that the  functions ${F_i}_{|U_i\cap H}$ are  constant whenever $U_i\cap H \ne \emptyset$.
Thus  if $i:H \to M$ denotes the inclusion then  $i^* \omega = 0$. This proves assertion (1).

Assertion (2) follows from similar considerations.
\end{proof}

\begin{definition}\label{D:FD}\rm
If  $D\neq 0$ is a $S^1$-flat divisor of a compact complex manifold $M$ then  $\mathcal F_D$ is the codimension one
singular holomorphic foliation induced by $\omega_D$.
\end{definition}

Note that the foliation  $\mathcal F_D$ is defined in an unambiguous way thanks to  item (2) of proposition \ref{P:basico}.
A key property of the foliation $\mathcal F_D$ is described in the following

\begin{cor}\label{C:folheacao}
If  $D \neq 0$ is a $S^1$-flat line bundle over a complex manifold $M$ then  the foliation $\mathcal F_D$
 leaves $D$ and  every compact complex hypersurface contained in the complement of the support of $D$ invariant.
\end{cor}

Let $M$ be the surface describe in example \ref{E:Hopf2} and  $D$ be  a non-zero  divisor on $M$ supported on the quotient
of one of the axis. Since the quotients of the two axis do not intersect each other it follows from corollary \ref{C:folheacao} that
$D$ is not $S^1$-flat, although  $D$ is $\mathbb C^*$-flat.

Another consequence of proposition \ref{P:basico} is a particular case of Theorem \ref{T:fibracao}. We will include it here since
the arguments are simpler then in the general case. 

\begin{cor}\label{C:simplificado}
Let $M$ be a compact complex manifold  such  that  $2h^1(M,\mathbb R)=h^1(M,\mathcal O_M)$.
If  $\{H_i\}_{i \in \mathbb N}$  is an infinite set of pairwise disjoint hypersurfaces of $M$ then there exists a
holomorphic map $\rho:M\to C$, $C$ a smooth algebraic curve, with connected fibers such that every $H_i$ is a component
of a fiber of $\rho$.
\end{cor}
\begin{proof}Since $H^2(M,\mathbb R)= H^2(M,\mathbb Q) \otimes \mathbb R$ is finite dimensional there exists  integers $k$, $n_1, \ldots, n_k$, $k>0$,
such that
\[
   c_{\mathbb R}\left(  \mathcal O_M(\sum_{r=1}^k n_r H_r )  \right) = 0 \, ,
\]
i.e., the divisor $D=\sum_{r=1}^k n_r H_r $ has zero real Chern
class. From proposition \ref{P:flat}(items (3) and (5)) we obtain
that the line bundle $\mathcal O_M(\sum_{i=1}^k n_r H_r )$ is
$S^1$-flat. From corollary \ref{C:folheacao} there exists a
codimension one holomorphic foliation $\mathcal F$ of $M$ leaving
every $H_i$, $i \in \mathbb N$, invariant.

We now make use of of Ghys' version of Jouanolou's Theorem, see \cite{Ghys}, to obtain a meromorphic first integral
$g: M \dashrightarrow P^1_{\mathbb C}$ of $\mathcal F$.
Since the hypersurfaces are $\{H_i\}$ are pairwise disjoint we can easily verify that the indeterminacy locus of
$g$ is empty and therefore $g$ is holomorphic.
>From Stein's factorization Theorem there exists an algebraic curve $C$, a fibration with
 $\rho: M \to C$ and a ramified covering $\pi:C\to P^1_{\mathbb C}$  such that $g= \pi  \circ \rho$.
\end{proof}

We end  this section with examples of $S^1$-flat divisors over projective manifolds which are not associated to
divisors linearly equivalent to zero.

\begin{example}\rm
Let $M$ be projective manifold and suppose that there exists a nontrivial homomorphism $\phi: \pi_1(M) \to S^1$.
If $\pi: \tilde M \to M$ is the universal covering of $M$  then we consider the codimension one foliation $\tilde{ \mathcal G}$ over $M \times \mathbb C^2$
defined by the $1$-form $\omega = xdy -ydx$ where $(x,y)$ are the coordinates of $\mathbb C^2$.
The homomorphism $\phi$ induces and action $ \Phi$ of $\pi_1(M)$ on $\tilde M \times \mathbb C^2$ given by
\begin{eqnarray*}
\Phi: \pi_1(M) \times (\tilde M \times \mathbb C^2) &\to& \tilde M \times \mathbb C^2  \\
(g,(p,(x,y))) &\mapsto& (g\cdot p, (\phi(g) x , \phi(g)^{-1} y))
\end{eqnarray*}
It is easy to see that the action $\Phi$ preserves the foliation $\tilde {\mathcal G}$.

The quotient of $\tilde M \times \mathbb C^2$ by $\Phi$ defines a rank $2$ vector bundle $E$ over $M$ equipped with a codimension one foliation $\mathcal G$.
Observe that $E= \mathcal L \oplus \mathcal L^*$ where $\mathcal L$ and $\mathcal L^*$ are flat line-bundles.

By GAGA's principle this vector bundle is algebraic  and therefore $\mathbb P (E)$, the projectivization of $E$, is a projective variety. Note that
$\mathbb P (E)$ carries two sections $M_1$ and $M_2$ corresponding to the split $E= \mathcal L \oplus \mathcal L^*$.

The foliation $\mathcal G$ induces a smooth codimension one foliation $\mathcal F$ of $\mathbb P (E)$ which leaves the two sections $M_1$ and $M_2$ of $\mathbb P (E)$
invariant. The divisor $D= M_1 - M_2$ is $S^1$-flat and it is possible to prove that $\mathcal F= \mathcal F_D$. Moreover we have that
\begin{enumerate}
\item if the image of $\pi_1(M)$ is a non-trivial  finite subgroup of $S^1$ then $D$ is not linearly equivalent to zero but there exist
a multiple of $D$ linearly equivalent to zero;
\item if the image of $\pi_1(M)$ is not
a finite subgroup of $S^1$ then $D$, or any of its multiples, is   not linearly equivalent to zero.
\end{enumerate}
\end{example}

The next example is a variant of an example presented  in \cite{Totaro} and attributed to Brendan Hassett.

\begin{example}\rm
Let $M$ be a projective variety with $h^1(M,\mathcal O_M) > 0$ and $\mathcal L$ be a non-trivial line-bundle with trivial Chern class, i.e., $\mathcal L \in Pic_0(M)$.
Let  $D_1$ be  an effective  divisor such that  $h^0(M,\mathcal O_M(D_1) \otimes \mathcal L)>0$. If $D_2$ is the zero divisor
of a section of  $\mathcal O_M(D_1) \otimes \mathcal L$ then $D = D_1 - D_2$ is $S^1$-flat. Moreover
\begin{enumerate}
\item if $\mathcal L$ is a torsion element of $Pic_0(M)$ then   $D$ is not linearly equivalent to zero but there exist
a multiple of $D$ linearly equivalent to zero;
\item if $\mathcal L$ is a non-torsion element of $Pic_0(M)$ then  $D$, or any of its multiples, is   not linearly equivalent to zero.
\end{enumerate}
\end{example}

\section{The Group $\Gamma(M)$ for compact complex manifolds}\label{S:Gamma}

If $M$ is a compact complex manifold let $\Gamma(M)$ be the group defined by
\[
\Gamma(M) = \frac{{\rm Div(M)}\otimes \mathbb Q}{{\rm Div}_{S^1}(M) \otimes \mathbb Q }.
\]
Since flat divisors have zero rational Chern class then there exists a rational Chern class map
\[
  c_{\mathbb Q}: \Gamma(M) \to H^2(M,\mathbb Q) \, .
\]
It follows from proposition \ref{P:flat} that for compact
complex manifolds $M$ with $h^1(M,\mathcal O_M)= 2h^2(M,\mathbb R)$  that the map $c_{\mathbb Q}$ is
injective. For projective manifolds  it is a trivial matter to verify that this injectivity identifies the
image of $\Gamma(M)$ with the rational Neron-Severi group of $M$.

In order to prove that $\Gamma(M)$ is finite diminensional for general compact complex manifolds we will
consider the algebraic reduction of  $M$.

\subsection{Divisors and Algebraic Reduction}\label{S:red}

For the results mentioned on the next two paragraphs the reader should consult \cite[pages 24--27]{Ueno} and references there within.

If $M$ is a compact complex manifold then the field of meromorphic functions of $M$, denoted by $k(M)$,
is a finitely generated extension of $\mathbb C$ whose transcendence degree is bounded by the dimension of $M$.
The transcendence degree of $M$ is called its algebraic dimension and will be denoted by $a(M)$. In the case
$a(M)= \dim M$ then $M$ is called a Moishezon manifold and there exists a finite succession of blow-ups along
non-singular centers such that the resulting manifold is projective.

In general there exists a compact complex variety $\tilde M$,  a   bimeromorphic morphism   $\psi: \tilde M \to M$
and  a   morphism $\pi: \tilde M \to N$ with connected fibers such that $N$ is a smooth projective  variety and
\[
  \psi^* k(M) = \pi^* k(N) \, .
\]
The projective variety $N$ is called an {\it algebraic reduction} of $M$. Note that
an algebraic reduction of $M$ is unique up to bimeromorphic equivalence.

We will say that a hypersurface $H$ of a complex variety $M$ is {\it special} if, in the notations above, the restriction
of $\pi \circ \psi^{-1}$ to $H$ is a dominant meromorphic map, i.e., has dense image. Remark that a Moishezon variety does not have special hypersurfaces
and every hypersurface of a variety of zero algebraic dimension is special.

The proposition below is a  generalization of Theorem 6.2 of \cite[page 129]{BPV}.

\begin{prop}\label{P:BPV}
If  $M$ is  a compact complex variety then there are at most $h^1(M,\Omega^1_{M})  + \dim M -  a(M) $ special hypersurfaces.
\end{prop}
\begin{proof}
Suppose that $H_k$, $1\le k \le   h^1(M,\Omega^1_{M})  + \dim M -  a(M) + 1$, are  distinct special  hypersurfaces of $M$
and let $\mathbb H = \oplus_{k}  \mathbb C \cdot H_k$ be the $\mathbb C$-vector space generated by them.

As in \cite{Ghys}
we can define a morphism from $\theta:\mathbb H \to H^1(M,\Omega^1_M)$ as follows: for every $H_k$ we can consider the associated
line bundle $\mathcal O_M(H_k)$ and map it to $d \log \phi_{ij}$, where $\phi_{ij}$ are the transition functions of $\mathcal O_M(H_k)$;
the morphism is them defined through linearity. If $\sum \lambda_k H_k$ belongs to the kernel of $\theta$ then we can define
a global meromorphic $1$-form with simple poles of residue $\lambda_k$ along $H_k$.

>From our assumptions we have that the dimension of the kernel of $\theta$ is at least $\dim M -  a(M) + 1$
and we can therefore construct $\omega_1, \ldots, \omega_l$,  $l=\dim M -  a(M) + 1$,
meromorphic $1$-forms over $M$ such that the polar set of $\omega_r$ is
\[
   (\omega_{r})_{\infty} = H_1 \cup H_2 \cup \ldots \cup H_{h} \cup H_{h+r} \, ,
\]
where $h=h^1(M,\Omega^1_{M})$. In particular the $1$-forms $\omega_i$ are linearly independent over $\mathbb C$ and moreover
since $H_i$ are special theirs restriction to $F$, the closure of a general fiber of  $\pi \circ \psi^{-1}:M\dashrightarrow N$,
are still  linearly independent over $\mathbb C$.

Let $\eta_1, \ldots, \eta_{a(M)}$ be the pullback under $\pi \circ \psi^{-1}$ of rational $1$-forms of $N$
linearly independent over $k(N)$. Since we have now $\dim M + 1$ meromorphic $1$-forms there exists
a relation of the form
\[
  \sum_{i=1}^l f_i \omega_i = \sum_{j=1}^{a(M)} g_j \eta_j \, ,
\]
where $f_i$ and $g_j$ are meromorphic functions of $M$.

If we take the restriction of the relation above to $F$ then since every meromorphic function of $M$
is constant along $F$   we obtain that the restriction of the $1$-forms  $\omega_i$ to $F$  are linearly
dependent over $\mathbb C$. A contradiction which proves that there are at most $h^1(M,\Omega^1_{M})  + \dim M -  a(M)$
special hypersurfaces on $M$.
\end{proof}

The statement  of Theorem 6.2 of \cite[page 129]{BPV} is the specialization to the case of surfaces of the following

\begin{cor}\label{C:BPV}
If  $M$ is  a compact complex variety of algebraic dimension zero
then $M$ has  at most $h^1(M,\Omega^1_{M})  + \dim M$  hypersurfaces.
\end{cor}

\begin{example}\rm
If $M$ is the quotient of $\mathbb C^n \setminus \{0 \}$ by a sufficiently general contraction then
$a(M)=h^1(M,\Omega^1_M)=0$ and $M$ has $n=\dim M$ special hypersurfaces; this shows that  the bound presented in proposition \ref{P:BPV}
above is sharp in every dimension.
\end{example}

\subsection{ Proof Theorem \ref{T:finito}}
If the algebraic dimension of $M$ is zero then Theorem \ref{T:finito} is an   immediate consequences of  corollary \ref{C:BPV}.

If $M$ is a Mosheizon variety then, as we have already mentioned, there exists
a projective variety $\tilde M$ and a bimeromorphic morphism $\psi: \tilde M \to M$.
Therefore  the finitude of $\dim_{\mathbb Q} \Gamma(\tilde M)$ implies the finitude of $\dim_{\mathbb Q} \Gamma( M)$.

It remains to deal with the cases where  the algebraic dimension of
$M$ satisfies $0 < a(M) < \dim M$.

Without loss of generality we can suppose that there exists a holomorphic map  from $M$
to a smooth projective variety $N$ with connected fibers. Let $\pi: M \to N$ be such map and set   $\mathcal R\subset M$ as
\[
  \mathcal R = \{ x \in M; {\rm rank} \, d\pi(x) < \dim N \} \, .
\]

We will make use of the following  lemma.

\begin{lemma}\label{L:elementar}
In the notations above if  $H$ be an irreducible hypersurface of $M$ then
\begin{enumerate}
\item if $\dim \pi(H) < \dim N -1$ then $H \subset \mathcal R$;
\item if $\pi(H) \not\subset \pi(\mathcal R)$ and $\dim \pi(H)=\dim N -1$  then, in the group of divisors of $M$, $\pi^{*} (\pi(H))=H + E$,
where $E$ is an effective divisor supported on $\mathcal R$.
\end{enumerate}
\end{lemma}
\begin{proof}
Item (1) follows from the local form of submersions and item (2) follows from Sard's Theorem and
 the connectedness  of the fibers of $\pi$.
We leave the details to the reader.
\end{proof}

Let now $S(M)$ denote the subgroup of ${\rm Div}(M)$ generated by the special divisors of $M$ and $R(M)$ denote the subgroup of ${\rm Div(M)}$  generated
by divisors supported on $\mathcal R$. Note that  $S(M)$ and $R(M)$ have both of finite rank and
the map
\begin{eqnarray*}
  {\rm Div}(N) \oplus S(M) \oplus R(M)  &\to& {\rm Div}(M) \,  \\
    (D,S,R) &\mapsto& \pi^* D + S + R
\end{eqnarray*}
is surjective. Since $\pi^*$ sends $S^1$-flat divisors of $N$ to
$S^1$-flat divisors of $M$ and  $\Gamma(N)$ is finite dimensional
then Theorem \ref{T:finito} follows. \qed

\section{A key property of the foliation $\mathcal F_D$ }

We start this section with  a simple lemma.

\begin{lemma}\label{L:Hironaka}
If $D$ is a divisor of a compact complex manifold $M$ then there exists a compact complex manifold
$\tilde M$ and a bimeromorphic map $\pi: \tilde M \to M$ such that $\pi^*D$ admits a decomposition of
the form $\pi^*D= D_{+} - D_{-}$ where $D_+$ and $D_-$ are effective divisors with disjoint supports.
\end{lemma}
\begin{proof}
>From  Hironaka's desingularization Theorem we can suppose that the support of $D$ has only normal crossing singularities
and its irreducible components are smooth.

Write $D$ as $P_0 - N_0$  where $P_0$ and $N_0$ are effective divisors without irreducible components
in common in theirs supports. If $\mathcal V_0$ be the set of codimension two subvarieties $V$ of $M$ such that
$V \subset H_P \cap H_N$ where $H_P$ is an irreducible component of $P_0$ and $H_N$ is an irreducible component of $N_0$.
Since we are supposing that the support of $D$ has normal crossings then every $V \in \mathcal V_0$ is smooth.

The proof follows from an induction on the cardinality of $\mathcal V_0$. We leave
the details to the reader.
\end{proof}

>From now on we will say that {\it a  divisor $D$ is without base points} if $D=D_+ - D_-$ where
 $D_+$ and $D_-$ are effective divisors with disjoint supports.

The next  proposition  is  the cornerstone of the proofs of Theorems
\ref{T:fibracao} and \ref{T:FolhasGeral}. It  is in fact a generalization of item (1) of
proposition \ref{P:basico}.

\begin{prop}\label{P:final}
Let $\mathcal F$ be a holomorphic foliation of a compact complex manifold $M$ and $D$ be a $S^1$-flat divisor of $M$.
If $\mathcal F$ admits a transcendental leaf which do not intersect the support of $D$ then $\mathcal F = \mathcal F_D$.
\end{prop}
\begin{proof}
>From lemma \ref{L:Hironaka} we can suppose, without loss of generality, that the $S^1$-flat divisor $D$ is without base points, i.e.,
$D=D_+ - D_-$ where $D_+$ and $D_-$ are effective divisors with disjoint supports.

By definition  the foliation $\mathcal F_D$ is induced by a closed meromorphic one form $\omega_D$
with polar set supported on $D$ and admits a real first integral $F:M \to [0,\infty]$ such that $F^{-1}(0)$ is equal to the
support of $D_+$ and $F^{-1}(\infty)$ is equal to the
support of $D_-$, see the proof of proposition \ref{P:basico} and definition \ref{D:FD}.

Consider
\[
  F_{|L}: L \to [0 , \infty ]
\]
the restriction of $F$ to the transcendental leaf $L$ of $\mathcal F$. Since $F$ is locally defined as
the modulus of a holomorphic function then $F_{|L}$ is either constant or an open map.
If $F_{|L}$ is constant then $L$ is a leaf of both $\mathcal F$ and $\mathcal F_D$. In particular the tangency locus
of $\mathcal F$ and $\mathcal F_D$ contains the analytic closure $L$. Since $L$ is not contained in any hypersurface
then $\mathcal F= \mathcal F_D$.

 From now on we will suppose that $F_{|L}$ is an open map. Recall that $F$ can be locally written as
\[
   F = \left| \exp \left( \int \omega_D \right) \right|
\]
and that all the periods of $\omega_D$ are purely imaginary complex numbers. If all the periods of $\omega_D$
are commensurable with $\pi \sqrt{-1}$ then there exists a positive integer $n$ such that $F^n$ is equal
to the modulus of the complex function $f=\exp( n ( \int \omega_D ))$.  Let $\overline L$ be the topological
closure of $L$ and $\partial L= \overline L \setminus L$.  Since $f$ is open and $\overline L$ does not intersect
the support of $D$  then $f(\partial L) = \overline{ f(L)} \setminus f(L)$. Therefore  $\partial L$ is invariant by both $\mathcal F$
and $\mathcal F_D$. Since $f(L)$ is an open set  relatively compact in $\mathbb C^* \subset \mathbb P^1$ we have that
$\partial f(L)$ is infinite. In particular $\mathcal F$ and $\mathcal F_D$ have an infinite number of leaves in common. This is sufficient
to show that $\mathcal F = \mathcal F_D$.

It remains to analyze the case where  $\omega_D$ has at least two linearly independent periods.  When this is the case the
multi-valued function $f=\exp(  \int \omega_D )$ has a monodromy group dense in $S^1$. Let $\pi: \tilde{M} \to M$ the
covering of $M$ associated to $f$ and consider the commutative diagram below.
\[
\begin{CD}
\pi^{-1} ( L) @>>>  \tilde M   @>f>>   \mathbb P^1  \\
@V\pi VV  @V\pi VV                                    @VV|\cdot|V \\
 L @>>>  M                   @>F>>      [0,\infty]
\end{CD}
\]
Since $F_{|L}$  is open and $\overline L$ does not intersect the support of $D$ then there exist positive real numbers  $N_-$ and $N_+$ such that
$F(L) = (N_-,N_+)$. The density of the monodromy group of $f$ in $S^1$ implies that
\[
   f(\pi^{-1}(L)) =  \{ z \in \mathbb C ; N_-< |z| < N_+ \} .
\]
It follows easily that $\partial L$ contains an infinite number of leaves of $\mathcal F$. This is sufficient to show that $\mathcal F= \mathcal F_D$.

\end{proof}

\section{Compact Complex Varieties which fiber over a curve}\label{S:compacta}

The main purpose of this section is to prove Theorem \ref{T:fibracao}. The proof will be based on
the following lemma which is a corollary to proposition \ref{P:final}.

\begin{lemma}\label{L:fibracao}
Let $D_1 $ and $D_2$ be two $S^1$-flat divisors on a compact
complex manifold $M$. If there exists a connected component of the
support of $D_1$ which does not intersect the support of $D_2$
then $\mathcal F_{D_1}= \mathcal F_{D_2}$.
\end{lemma}
\begin{proof}
Let $\omega_{D_1}$ and $\omega_{D_2}$ be the meromorphic $1$-forms canonically associated to $D_1$ and $D_2$ and let $V$ be the connected component
of the support of $D_1$ which does not intersect the support of $D_2$.
In the proof of proposition \ref{P:final} we saw that  $\mathcal F_{D_1}$ admits a real analytic first integral $F:M\to [0,\infty]$ and $V$ is a
level of $F$.  Thus  there exists an open neighborhood $U$ of $V$ in $M$ which is
saturated by the foliation $\mathcal F_{D_1}$, i.e., every leave of $\mathcal F_{D_1}$ which intersects $U$ is in fact contained in $U$. Moreover
we can choose $U$ such that $D_2 \cap U = \emptyset$.

Let $L$ be an arbitrary leaf of $\mathcal F_{D_1}$ contained in $U \setminus V$.
If for every such $L$, $\overline L$ is a complex subvariety of $M$ then
$\overline L$ is invariant by  $\mathcal F_{D_2}$ by item (1) of proposition \ref{P:basico}. Thus the restriction of $\mathcal F_{D_1}$ and $\mathcal F_{D_2}$
to $U$ coincide and therefore  $\mathcal F_{D_1}= \mathcal F_{D_2}$.
Otherwise  $F_{D_1}$ has a transcendental leaf which does not intersect the support of $D_2$ and  the lemma follows from proposition \ref{P:final}.
\end{proof}

\subsection{A foliation theoretic proof  of Totaro's  Theorem}\label{S:Totaro}
>From lemma \ref{L:fibracao} we obtain a slight generalization of part
of a result by Totaro \cite[Theorem 2.1]{Totaro}. The original result
is stated for complex projective manifolds and its original proof has the advantage to work over
fields of positive characteristic, see \cite[section 6]{Totaro}. Here we will adopt a foliation theoretic approach
which works uniformly for every compact complex variety(and give no hint how to proceed to establish the Theorem in positive characteristic).

\subsubsection{\rm \bf Proof of Theorem \ref{T:Totaro2}} From the hypothesis of the Theorem
there exists integers $n_{12}, n_{23}, m_{12}, m_{23}$ such that
\[
 D_{12}=n_{12} D_1 - m_{12} D_2 \, \, \text{and} \, \, D_{23}=n_{23} D_2 - m_{23}
D_3
\]
are $S^1$-flat divisors without base points. Lemma
\ref{L:fibracao} implies that  $\mathcal F_{D_{12}} = \mathcal
F_{D_{23}}$.

We can now conclude as in the proof of Jouanolou's Theorem. Note
that there exists a non-constant meromorphic function $F \in k(M)$
such that $\omega_{D_{12}}=F \cdot \omega_{D_{23}}$.
Differentiating we obtain that $dF \wedge \omega_{D_{23}} = 0$,
i.e., $F$ is a meromorphic first integral of $\mathcal
F_{D_{12}}=\mathcal F_{D_{23}}$. Since $D_{12}$ is without base
points it follows that $F$ is in fact a holomorphic map $F:M \to
\mathbb P^1$. The Theorem follows taking the Stein factorization
of $F$.
 \qed

\subsection{Proof of Theorem \ref{T:fibracao}}
Considering the natural map $\mathcal H\otimes \mathbb Q \to
\Gamma(M)$ it follows that there exists two $S^1$-flat divisors
$D$ and $D'$ with support contained in $\mathcal H$ and such that
there exists a component $H$ of the support of $D$ not contained
in $D'$. At this point we can use lemma \ref{L:fibracao} to
guarantee that $\mathcal F_{D} = \mathcal F_{D'}$ and conclude as
in Jouanolou's Theorem, cf. the proof of Theorem \ref{T:Totaro2}.
 \qed

\section{Diffeomorphism type  of Smooth Divisors}
This section is devoted to the proof of Theorem \ref{T:Totaro}.
Before proceeding to the proof we would like to recall some remarks and examples made by Totaro in \cite{Totaro}.
\begin{enumerate}
\item If $M$ is a projective variety with $H^1(M,\mathbb R)=0$ then the Betti numbers and Hodge numbers of a smooth divisor
are  determined by its Chern class, see \cite[remark 2]{Totaro}. More generally the same argument used there
apply to any  compact complex variety with $H^1(M,\mathcal O_M)=0$.
\item If $M$ is a projective and $D$ is an ample smooth divisor of $M$ then the Betti and Hodge numbers of $D$ are determined by
its Chern class, see \cite[remark 1]{Totaro}.
\item There exists smooth complex projective manifolds with two disjoint homologous smooth divisors which are both connected
but have different Betti numbers.
\end{enumerate}

\subsection{Proof of Theorem \ref{T:Totaro}} Let $D_1$ and $D_2$ be two connected divisors whose $\Gamma$-classes lie in a line of
$\Gamma(M)$. There exists integers $p$ and $q$ such that $D = pD_1
- qD_2$ is $S^1$-flat. Thus we can  choose a covering $\mathcal U=
\{U_i\}$ os $M$ and local rational functions $f_i$ defining
$D=pD_1 - qD_2$ such that the $f_i= f_{ij} f_j$ and $f_{ij}$ are
locally constant functions of modulus $1$. As before we have that
$\omega = \frac{df_i}{f_i}$ is a well-defined global meromorphic
$1$-form and moreover $F=|f_i|: M \to [0, \infty]$  is a
well-defined continuous function(real-analytic outside the support
of $D$).

We can also define a global (real) $1$-form $\theta$
over $M\setminus (D_1 \cup D_2)$ by the relation $\theta_{|U_i} = d \arg (f_i)$, where $\arg$ denotes the complex argument, i.e.,
$f = |f|\cdot \exp(\arg(f))$. In a neighborhood of  $D_1$ and $D_2$ the $1$-form $\theta$ has mild algebraic singularities; if $U$ is a
sufficiently small neighborhood of a point $p\in D_1 \cup D_2$   then, over $U$, the foliation induced by $\theta$ is diffeomorphically equivalent
to  the foliation of $\Sigma\times (D\cap U)$ induced by $xdy - ydx$, where $(x,y)$ are local real coordinates of a transversal $\Sigma$  of $D$.

Integration along closed paths defines a homomorphism
\begin{eqnarray*}
   \int \theta: H_1(M\setminus (D_1 \cup D_2),\mathbb R) &\to& \mathbb R  \, \\
   \gamma &\mapsto& \int_{\gamma} \theta
\end{eqnarray*}
which sends $\gamma_1$ and $\gamma_2$, small loops around $D_1$ and $D_2$ respectively, to real numbers commensurable to
$\pi$.

The inclusion of $M\setminus(D_1 \cup D_2)$ into $M$ induces a surjective  homomorphism
$H_1(M \setminus(D_1 \cup D_2) , \mathbb R)  \to H_1(M  , \mathbb R)  $    whose kernel is contained in the subspace of $H_1(M \setminus(D_1 \cup D_2) , \mathbb R)$
generated by $\gamma_1$ and $\gamma_2$.
Therefore we can choose a morphism $T:H_1(M  , \mathbb R) \to \mathbb R$ such that
$(T + \int \theta)(\gamma)$ is commensurable to $\pi$ for every $\gamma \in  H_1(M \setminus(D_1 \cup D_2) , \mathbb Z)$. From DeRham's
isomorphism we deduce the existence of a (real) closed $1$-form $\eta$ on $M$  such that $\int_{\gamma} \theta + \eta $ is a rational multiple
of $\pi$ for every $\gamma \in  H_1(M \setminus(D_1 \cup D_2) , \mathbb Z)$.
Since $H_1(M \setminus(D_1 \cup D_2) , \mathbb Z)$ is finitely generated there exists an integer $N$ such that
\[
     G=\exp\left (  i N \int \theta + \eta \right): M \setminus(D_1 \cup D_2)\to S^1
\]
 is a well-defined $C^{\infty}$ function.

Since the local structure, around points of $D_1$, of the
foliation induced by $\eta + \theta$ is the same of the foliation
induced by $\theta$ if we take a smooth fiber $\tilde{D_1}$ of
\[
   F^N\times G : M \setminus(D_1 \cup D_2) \to (0,\infty) \times
S^1 \cong \mathbb C^*
\]
 sufficiently close to $D_1$ then  there exists a positive
integer $m$ such that $\tilde{D_1}$ is an \'etale $\mathbb Z /
m$-covering of $D_1$. In an analogous way a smooth fiber
$\tilde{D_2}$ of $F\times G$ sufficiently close to $D_2$ is an
\'etale $\mathbb Z / n$-covering of $D_2$ for some positive
integer $n$.

We claim that  choosing $\eta$ small enough we can join
$q_1=(F^N\times G)(\tilde D_1)$ to $q_2=(F^N\times G)(\tilde D_2)$
by a differentiable path $\gamma:[ 0,1] \to \mathbb C^*$ avoiding
the critical values of $F^N\times G$. In fact if
\[
\rho : \overline M \to M \setminus ( D_1 \cup D_2 ) \,
\]
 denotes the universal covering of $M \setminus ( D_1 \cup D_2 )$
 then $\rho^* \theta$ is exact and holomorphic. Thus $\rho^*(N \omega)$ admits a holomorphic primitive
 \[
   \overline G: \overline M \to \mathbb C
\]
which has a set of critical values of real codimension two. Thus
if $\eta$ is small enough then $\rho^* (F^N\times G) $ is
sufficiently close to $\overline G$ and thus the set of critical
values are also close. This is sufficient to prove the claim.

We can therefore lift the real vector field
$\frac{\partial}{\partial t}$ on $[0,1]$ to $(F^N\times
G)^{-1}(\gamma([0,1]))$. The flow of this vector field induces  a
diffeomorphism between $\tilde D_1$ and $\tilde D_2$.  \qed

\begin{remark} \rm
 Totaro  shows in \cite{Totaro} that if $D_1$ and $D_2$ are smooth
divisors on a projective manifold with the same Chern class then,
after blowing up the intersection scheme of $D_1$ with $D_2$ and
resolving the resulting variety, the strict transforms of $D_1$
and $D_2$ have the same Chern class. A similar argument shows that
our result holds for smooth divisors $D_1$ and $D_2$, not
necessarily disjoint, with the {\bf same} $\Gamma$-class.

Note that a similar reasoning for smooth divisors which are not
disjoint and whose $\Gamma$-classes lie in a line does not work.
For instance if we take the a line $L$ and a smooth cubic $C$ on
$\mathbb P^2$ they have non-diffeomorphic universal coverings and
 $3L-C$ is a $S^1$-flat divisor.
\end{remark}

\section{The Closure of Transcendental Leaves}\label{S:Folhas}

\subsection{Hypersurfaces with ample normal bundle}

The original motivation of this work was to find an analogous of
the following well-known fact for general projective varieties: if
$\mathcal F$ is a holomorphic foliation of $\mathbb P ^n$ and $L$
is a transcendental leaf of $\mathcal F$ it is well known that the
topological closure of $L$ does intersect every compact
hypersurface of $\mathbb P^n$, see \cite{GB}.

The key point on the proof of the fact above is that the
complement of any compact hypersurface of $\mathbb P^n$ is Stein,
and even more it is in fact affine.

A first result, and almost obvious, result on this direction is
the following

\begin{prop}\label{P:inter}
Let $\mathcal{F}$ be a holomorphic foliation of a projective variety $M$
and  $H$ an effective divisor  of $M$ with ample normal bundle.
\begin{footnote}{See \cite{Hart} for a precise definition of divisors with ample normal bundle.}
\end{footnote}
 If $L$ is a leaf
of $\mathcal{F}$  then   $L$ is contained in contractible subvariety
of $M$ or the topological closure of $L$  intersects $H$.
\end{prop}

\begin{proof}
Since the normal  bundle of $H$ is ample according to  \cite[Theorem 4.2, p. 110]{Hart}
there exists a positive integer
$k$ such that:

\begin{itemize}
\item the linear system ${\rm H}^0(M,\mathcal{O}_M(nH))$
is free from base points;
\item  the natural map $\phi:M \to {\mathbb
P}{\rm H}^0(M,\mathcal{O}_M(nH))$ is holomorphic
\item $\phi$ is  biholomorphic
on a neighborhood of $H$;
\item the set $\phi(H)$  can be identified with the
intersection of an hyperplane of ${\mathbb P}{\rm
H}^0(M,\mathcal{O}_M(nH))$ with $\phi(M)$, the image of $\phi$.
\end{itemize}

Thus we can
identify $\phi(M \setminus H)$ with an
affine closed set of some affine space
$\mathbb{C}^N$.

Suppose that $L$ is not contracted by $\phi$.  Therefore we can choose  a principal
open subset $U$ of $\mathbb{C}^N$  satisfying:

\begin{itemize}
\item $U \cap \phi(L) \neq \emptyset$;
\item  $\phi(M \setminus H) \cap U$ is a smooth affine variety;
\item  $\phi_* \mathcal{F}$ restricted to $U \cap \phi(M \setminus H)$ is generated by
a global section of the tangent sheaf of $U\cap \phi(M \setminus H)$.
\end{itemize}

Denote by $\tilde{M}$ the intersection of $\phi(M \setminus H)$ with $U$.

Let $\Theta_{U}$ be the sheaf of vector fields of $U$ and $\Theta_{U,\tilde{M}}$ be the subsheaf formed
by vector fields tangent to $\tilde{M}$. From the exact sequence of $\mathcal{O}_U$-coherent sheaves
\[
  0 \to \Theta_{U,\tilde{M}} \otimes I_{\tilde{M}} \to \Theta_{U,\tilde{M}
  } \to \Theta_{\tilde{M}} \to 0 \, ,
\]
and the fact that coherent sheaves over affine varieties have no
higher order cohomology we deduce that there exists a foliation of $U$, which
naturally extends to a foliation $\mathcal{G}$ of
${\mathbb P}{\rm H}^0(V,\mathcal{O}_V(nH))$, whose restriction to
$\tilde{M}$ coincides with $\phi_*(\mathcal{F})$.

Therefore let $p \in \phi(L)\cap \tilde{M}$ be a non-singular point of $\phi_*(\mathcal{F})$.
 The
proposition  follows if the closure of the leaf of $\mathcal{G}$ through $p$ is not contained on any
compact subset of $\mathbb{C}^N$. But this is precisely the case
since no holomorphic vector field on $\mathbb{C}^N$ has a bounded
leaf.
\end{proof}

In the two dimensional case the above proposition specializes to

\begin{cor}\label{C:inter}
Let $\mathcal{F}$ be a holomorphic foliation of a projective surface  $S$
and  $D$ an effective divisor  of $S$. If $L$ is a transcendental leaf
of $\mathcal{F}$  such that topological closure of  $L$  does not intersect the support of  $D$ then
$D^2 \leq 0$.
\end{cor}

\subsection{Proof of Theorem \ref{T:FolhasGeral}}Let $\mathcal H$ be the set of compact complex hypersurfaces
of $M$ which do not intersect the leaves of $\mathcal F$. Denote
by $\Gamma$ the natural map
\[
  \Gamma : \mathcal H\otimes \mathbb Q \to \Gamma(M) \, .
\]

If $\dim_{\mathbb Q} \ker \Gamma \ge 2$ then it follows that there
exists two
 there exists two $S^1$-flat
divisors $D_1$ and $D_2$ such that the supports of $D_1$ and $D_2$
are  distinct. From proposition \ref{P:final} we have that
$\mathcal F=\mathcal F_{D_1} =\mathcal F_{D_2}$. Moreover, see
argument in the proof of the first part of Theorem
\ref{T:FolhasGeral}, $\mathcal F$ admits a meromorphic first
integral and do not admit transcendental leaves. This shows that
$\dim_{\mathbb Q} \ker \Gamma \le 1$. In  particular the
cardinality of $\mathcal H$ is at most $\dim_{\mathbb Q}
\Gamma(M)+1$.

When the cardinality of $\mathcal H$ is precisely $\dim_{\mathbb Q}
\Gamma(M)+1$ we have that $\dim_{\mathbb Q} \ker \Gamma=1$. If $D$
is a generator of the kernel of $\Gamma$ it follows from
proposition \ref{P:final}  that $\mathcal F=\mathcal F_{D}$. This
sufficient to prove (1).

If $M$ is projective then we claim that $\Gamma$ is not
surjective. Otherwise  there exists an ample  divisor $Z$ with
support contained in $\mathcal {H}$ contradicting  proposition
\ref{P:inter}. Moreover, since $M$ is projective,  if
$h^1(M,\mathcal O_M)=0$ then numerical and linear equivalence
coincides modulo torsion. Thus arguing as above we prove (2). \qed

We conclude by remarking that we cannot replace in the statement
of item (1) of Theorem \ref{T:FolhasGeral} the group $\Gamma(M)$
by the Neron-Severi group of $M$.

For instance if $M$ is an arbitrary  primary Hopf surface, i.e.
$S$ is the quotient of ${\mathbb C^*}^2$ by a linear diagonal
contraction, and  $E_1$ and $E_2$ denote the elliptic curves on
$S$ obtained as the quotients of the coordinates axis then the
divisor $D=E_1 - E_2$ is $S^1$-flat as the reader can easily
verify. If the algebraic dimension of $S$ is zero then every leaf
of the foliation $\mathcal F_D$ distinct from $E_1$ and $E_2$ is
transcendental and its  topological closure does not intersect the
support of $D$. Therefore for every transcendental leaf of
$\mathcal F_D$ the set $\mathcal H$ has cardinality two while the
Neron-Severi group has dimension zero, i.e.,  is the trivial
group.

\pagebreak
\begin{appendix}

\section*{Appendix by Laurent Meersseman}

It is a natural question to ask if the results of this paper are
still true if we replace $\Gamma (M)$ by the Neron-Severi group of
$M$ in the statements of the Theorems. The aim of this appendix is
to give a negative answer to this question, at least for Theorems
2,3 and 4. In other words, the group $\Gamma (M)$ is
really the good object to consider in these problems.

\begin{THM} There exists a compact, complex $5$-manifold $N$ with three
pairwise disjoint smooth hypersurfaces $H_1$, $H_2$ and $H_3$ such
that
\begin{enumerate}
\item[(a)] The Neron-Severi group of $N$ is reduced to zero.
\item[(b)] The manifold $N$ does not admit a holomorphic map onto
a smooth curve. \item[(c)] The universal coverings of $H_1$ and
$H_3$ are not homotopically equivalent.
\end{enumerate}
\end{THM}

The example of the Theorem comes from the family of compact,
complex manifolds constructed and studied in \cite{3} as a
generalization of \cite{2}. Let us first recall very briefly this
construction. Let $n>2m$ be positive integers. Let
$\Lambda=(\Lambda_1,\hdots,\Lambda_n)$ be a configuration of $n$
vectors of $\mathbb C^m$. Assume it is {\it admissible}, i.e. that it
satisfies

\noindent - the Siegel condition: $0\in\mathbb C^m$ belongs to the
(real) convex hull of $(\Lambda_1,\hdots,\Lambda_n)$.

\noindent - the weak hyperbolicity condition: if $0$ belongs to
the convex hull of a subset of $(\Lambda_1,\hdots,\Lambda_n)$,
then this subset has cardinal strictly greater than $2m$.

 Consider the holomorphic foliation $\mathcal F$ of the projective space $\mathbb P^{n-1}$
given by the following action
$$
(T,[z])\in\mathbb C^m\times\mathbb P^{n-1}\longmapsto [\exp\langle
\Lambda_1,T\rangle\cdot z_1,\hdots,\exp\langle
\Lambda_n,T\rangle\cdot z_n]\in\mathbb P^{n-1}
$$
where the brackets denote the homogeneous coordinates in $\mathbb
P^{n-1}$ and where $\langle -,- \rangle$ is the {\it inner}
product of $\mathbb C^n$. Define
$$
\mathcal N_{\Lambda}=\{[z]\in\mathbb P^{n-1}\quad\vert\quad
\sum_{i=1}^n\Lambda_i\vert z_i\vert ^2=0\}
$$
which is a smooth manifold due to the weak hyperbolicity
condition. Then, there exists an open dense subset $V\subset \mathbb
P^{n-1}$ such that the restriction of $\mathcal F$ to $V$ is
regular and admits $\mathcal N_{\Lambda}$ as a global smooth
transverse. Therefore, $\mathcal N_{\Lambda}$ can be endowed with
a structure of (compact) complex manifold as leaf space of
$\mathcal F$ restricted to $V$. We denote by $N_{\Lambda}$ this
compact complex manifold. It has dimension $n-m-1$ and is not
Kaehler if $n>2m+1$ (see \cite{3}, Theorem 2).

The standard action of the torus $(S^1)^n$ onto $\mathbb C^n$ leaves
$\mathcal N_{\Lambda}$ invariant and the corresponding quotient
space is easily seen to identify with a simple convex polytope
(see \cite{1}, Lemma 0.11; simple means dual to a simplicial
polytope). We denote by $P_{\Lambda}$ the combinatorial type of
this convex polytope. It has some remarkable properties:

\noindent (i) {\it Rigidity:} there is a $1:1$ correspondence
between the combinatorial classes of simple convex polytopes and
the classes of manifolds $N_{\Lambda}$ up to $C^{\infty}$
equivariant diffeomorphism and up to product by circles, see
\cite{1}, Theorem 4.1.

\noindent (ii) {\it Realization:} given any simple convex polytope
$P$, there exists $N_{\Lambda}$ such that $P_{\Lambda}=P$, see
\cite{3}, Theorem 13.

\noindent (iii) {\it Submanifolds:} a codimension $p$ face $F$ of
$P_{\Lambda}$ corresponds to a codimension $2p$ holomorphic
submanifold $N_{\Lambda'}$ of $N_{\Lambda}$ such that
$P_{\Lambda'}=F$, see \cite{3}, \S V.

We are now in position to construct our example.

\begin{proof} Consider the convex polyhedron $P$ obtained from
the cube by cutting off two adjacent vertices by a plane (cf
\cite{1}, Example 11.5). It has two triangular facets
(corresponding to the vertices which were cut off), two
rectangular ones, two pentagonal ones and finally two hexagonal
ones. By (ii), there exists manifolds $N_{\Lambda}$ such that
$P_{\Lambda}=P$. Here $N_{\Lambda}$ can moreover be assumed to be
$2$-connected (see \cite{3}, Theorem 13). Then, $n$ is equal to
$8$, and $m$ to $2$ so that $N_{\Lambda}$ has complex dimension
$5$. By (i), for all such choices of $\Lambda$, the
$C^{\infty}$-diffeomorphism type of $N_{\Lambda}$ is the same. Fix
such a $\Lambda$ and set $N=N_{\Lambda}$. By Corollary 4.5 of
\cite{1}, we may assume that $\Lambda$ is generic (in the sense of
condition (H) of \cite{3}, IV). The two triangular facets
correspond by (iii) to two smooth hypersurfaces $H_1$ and $H_2$.
Choose a rectangular facet of $P$, and let $H_3$ denote the
corresponding hypersurface of $N$. Notice that $H_1$, $H_2$ and
$H_3$ are pairwise disjoint since the corresponding facets are
pairwise disjoint.

Since $N$ is $2$-connected, its Neron-Severi group is reduced to
zero. This proves (a).

Since $\Lambda$ is generic, by \cite{3}, Corollary of Theorem 4,
then the manifold $N$ does not have any non-constant meromorphic
function. So cannot admit a holomorphic projection onto an
algebraic curve. This proves (b).

 Finally, using (i) and \cite{3}, VIII,
we have that $H_1$ and $H_2$ are diffeomorphic to $S^5\times
(S^1)^3$, whereas $H_3$ is diffeomorphic to $S^3\times S^3\times
(S^1)^2$. Now, the universal coverings of these two manifolds are
not homotopy equivalent. This finishes the proof. $\square$
\end{proof}

Notice that Theorem 4 implies that the rank of $\Gamma(N)$ is
greater than or equal to $2$.

\end{appendix}

\vskip.4cm

\noindent \small{Jorge Vit{\'o}rio Pereira\\
I.R.M.A.R.  - Universit\'e
de Rennes
I \\Campus de Beaulieu \\ 35042 Rennes Cedex  France \\
email: jvp@impa.br}

\end{document}